\documentclass[a4paper,12pt]{amsart}
\usepackage{amssymb}
\usepackage{ifthen}
\usepackage{graphicx}
\nonstopmode \numberwithin{equation}{section}
\usepackage[usenames]{color}
\newtheorem{thm}{Theorem}[section]
\newtheorem{cor}{Corollary}[section]
\newtheorem{lem}{Lemma}[section]

\newtheorem{rem}{Remark}[section]
\newtheorem{rems}[equation]{Remarks}
\theoremstyle{definition}
\newtheorem{defin}{Definition}[section]
\newtheorem{examp}[equation]{Example}
\newtheorem{prob}[equation]{Problem}
\newtheorem{ques}[equation]{Question}
\newtheorem{op}{Open Problem}[section]

\newtheorem{conj}[equation]{Conjecture}
\newtheorem{deter}[equation]{Determination}
\newtheorem{case}{Case}[section]
\newtheorem{subcase}[equation]{Subcase}
\newtheorem{claim}{Claim}[section]
\newtheorem{subclaim}{Subclaim}

\newcounter {own}
\def\theown {\thesection       .\arabic{own}}

\newenvironment{pf}[1][]{%
 \vskip 3mm
 \noindent
 \ifthenelse{\equal{#1}{}}%
  {{\bf Proof. }}%
  {{\bf #1.} }%
 }%
{\qed\bigskip}

\newcounter{alphabet}
\newcounter{tmp}
\newenvironment{Thm}[1][]{\refstepcounter{alphabet}%
\bigskip%
\noindent%
{\bf Theorem \Alph{alphabet}}%
\ifthenelse{\equal{#1}{}}{}{ (#1)}%
{\bf .} \itshape}{\vskip 8pt}

\makeatletter
\newcommand{\Ref}[1]{\@ifundefined{r@#1}{}{\setcounter{tmp}{\ref{#1}}\Alph{tmp}}}
\makeatother

\newcommand{\IR}{{\mathbb R}}

\newcommand{\IB}{{\mathbb B}}

\newcommand{\diam}{{\operatorname{diam}}}

\newcommand{\dist}{{\operatorname{dist}}}

\def\be{\begin{equation}}
\def\ee{\end{equation}}

\newcommand{\bee}{\begin{enumerate}}
\newcommand{\eee}{\end{enumerate}}

\newcommand{\blem}{\begin{lem}}
\newcommand{\elem}{\end{lem}}
\newcommand{\bthm}{\begin{thm}}
\newcommand{\ethm}{\end{thm}}
\newcommand{\bcor}{\begin{cor}}
\newcommand{\ecor}{\end{cor}}
\newcommand{\beg}{\begin{examp}}
\newcommand{\eeg}{\end{examp}}
\newcommand{\begs}{\begin{examples}}
\newcommand{\eegs}{\end{examples}}
\newcommand{\bdefe}{\begin{defin}}
\newcommand{\edefe}{\end{defin}}
\newcommand{\bprob}{\begin{prob}}
\newcommand{\eprob}{\end{prob}}
\newcommand{\bques}{\begin{ques}}
\newcommand{\eques}{\end{ques}}
\newcommand{\bei}{\begin{itemize}}
\newcommand{\eei}{\end{itemize}}

\newcommand{\bde}{\begin{deter}}
\newcommand{\ede}{\end{deter}}
\newcommand{\bca}{\begin{case}}
\newcommand{\eca}{\end{case}}
\newcommand{\bsca}{\begin{subcase}}
\newcommand{\esca}{\end{subcase}}
\newcommand{\bcl}{\begin{claim}}
\newcommand{\ecl}{\end{claim}}
\newcommand{\bscl}{\begin{subclaim}}
\newcommand{\escl}{\end{subclaim}}
\newcommand{\bcon}{\begin{conj}}
\newcommand{\econ}{\end{conj}}
\newcommand{\bcons}{\begin{conjs}}
\newcommand{\econs}{\end{conjs}}
\newcommand{\bprop}{\begin{propo}}
\newcommand{\eprop}{\end{propo}}
\newcommand{\br}{\begin{rem}}
\newcommand{\er}{\end{rem}}
\newcommand{\brs}{\begin{rems}}
\newcommand{\ers}{\end{rems}}
\newcommand{\bo}{\begin{obser}}
\newcommand{\eo}{\end{obser}}
\newcommand{\bos}{\begin{obsers}}
\newcommand{\eos}{\end{obsers}}
\newcommand{\bpf}{\begin{pf}}
\newcommand{\epf}{\end{pf}}
\newcommand{\ba}{\begin{array}}
\newcommand{\ea}{\end{array}}
\newcommand{\beq}{\begin{eqnarray}}
\newcommand{\beqq}{\begin{eqnarray*}}
\newcommand{\eeq}{\end{eqnarray}}
\newcommand{\eeqq}{\end{eqnarray*}}

\newcommand{\ds}{\displaystyle}

\newcommand{\bop}{\begin{op}}
\newcommand{\eop}{\end{op}}

\newtheorem{pfofThm1.5}[equation]{}

\newcounter{minutes}
\divide\time by 60
\newcounter{hours}
\multiply\time by 60 \addtocounter{minutes}{-\time}

\begin{document}

\bibliographystyle{amsplain}

\title{On bilipschitz extensions in real Banach spaces}

\def\thefootnote{}
\footnotetext{ \texttt{\tiny File:~\jobname .tex,
          printed: \number\year-\number\month-\number\day,
          \thehours.\ifnum\theminutes<10{0}\fi\theminutes}
} \makeatletter\def\thefootnote{\@arabic\c@footnote}\makeatother

\author{M. Huang $^* $}
\address{M. Huang, Department of Mathematics,
Hunan Normal University, Changsha,  Hunan 410081, People's Republic
of China} \email{mzhuang79@yahoo.com.cn}

\author{Y. Li}
\address{Y. Li, Department of Mathematics,
Hunan Normal University, Changsha,  Hunan 410081, People's Republic
of China} \email{yaxiangli@163.com}

\date{}
 \keywords{
Uniform domain, QH homeomorphism,
 bilipschitz map, bilipschitz extension.\\
${}^{\mathbf{*}}$ Corresponding author}

\begin{abstract}
Suppose that $E$ and $E'$ denote real Banach spaces with dimension
at least $2$, that $D\not=E$ and $D'\not=E'$ are bounded domains
with connected boundaries, that $f: D\to D'$ is an $M$-QH
homeomorphism, and that $D'$ is uniform.
 The main aim of this paper is to prove that $f$ extends to a homeomorphism $\overline{f}:
\overline{D}\to \overline{D}'$ and $\overline{f}|\partial D$ is
bilipschitz if and only if $f$ is bilipschitz in $\overline{D}$. The
answer to some open problem of V\"ais\"al\"a is affirmative under an
natural additional condition.
\end{abstract}

\thanks{The research was partly supported by NSFs of
China (No. 11071063 and No. 11101138).}

\maketitle{} \pagestyle{myheadings} \markboth{}{On bilipschitz
extensions in real Banach spaces}

\section{Introduction and main results}\label{sec-1}

During the past three decades, the quasihyperbolic metric has become
an important tool in geometric function theory and in its
generalizations to metric spaces and Banach spaces \cite{Vai5}. Yet,
some basic questions of the quasihyperbolic geometry in Banach
spaces are open. For instance, only recently the convexity of
quasihyperbolic balls has been studied in \cite{krt,rt} in the setup
of Banach spaces.

Our study is motivated by V\"ais\"al\"a's theory of freely
quasiconformal maps  and other related maps in the setup of Banach
spaces \cite{Vai6-0, Vai6, Vai5}. Our goal is to study some of the
open problems formulated by him. We begin with some basic
definitions and the statements of our results. The proofs and
necessary supplementary notation terminology will be given
thereafter.

Throughout the paper, we always assume that $E$ and $E'$ denote real
Banach spaces with dimension at least $2$. The norm of a vector $z$
in $E$ is written as $|z|$, and for every pair of points $z_1$,
$z_2$ in $E$, the distance between them is denoted by $|z_1-z_2|$,
the closed line segment with endpoints $z_1$ and $z_2$ by $[z_1,
z_2]$. We begin with the following concepts following closely the
notation and terminology of \cite{TV, Vai2, Vai, Vai6-0, Vai6} or
\cite{Martio-80}.

We first recall some definitions.

\bdefe \label{def1.3} A domain $D$ in $E$ is called $c$-{\it
uniform} in the norm metric provided there exists a constant $c$
with the property that each pair of points $z_{1},z_{2}$ in $D$ can
be joined by a rectifiable arc $\alpha$ in $ D$ satisfying

 \bee
\item\label{wx-4} $\ds\min_{j=1,2}\ell (\alpha [z_j, z])\leq c\,d_{D}(z)
$ for all $z\in \alpha$, and

\item\label{wx-5} $\ell(\alpha)\leq c\,|z_{1}-z_{2}|$,
\eee

\noindent where $\ell(\alpha)$ denotes the length of $\alpha$,
$\alpha[z_{j},z]$ the part of $\alpha$ between $z_{j}$ and $z$, and
$d_D(z)$ the distance from $z$ to the boundary $\partial D$ of $D$.
\edefe

\bdefe \label{def1.7} Suppose $G\varsubsetneq E\,,$ $G'\varsubsetneq
E'\,,$ and $M \ge 1\,$. We say that a homeomorphism $f: G\to G'$ is
{\it $M$-bilipschitz} if
 $$\frac{1}{M}|x-y|\leq |f(x)-f(y)|\leq M|x-y|$$
for all $x$, $y\in G$,  and {\it $M$-QH} if
 $$\frac{1}{M}k_{G}(x,y)\leq k_{G'}(f(x),f(y))\leq Mk_{G}(x,y)$$
for all $x$, $y\in G$. \edefe

As for the extension of bilipschitz maps in ${\mathbb R}^2$, Ahlfors
\cite{Ahl} proved that if a planar curve through $\infty$ admits a
quasiconformal reflection, it also admits a bilipschitz reflection.
Furthermore, Gehring gave generalizations of Ahlfors' result in the
plane.

\begin{Thm}\label{xwm-8}$($\cite[Theorem 7]{Ger6}$)$
Suppose that $D$ is a $K$-quasidisk in ${\mathbb R}^2$, that $D'$ is
a Jordan domain in ${\mathbb R}^2$ and that $\phi: \partial D\to
\partial D'$ is $L_1$-bilipschitz. Then there exist
$L$-bilipschitz $f: \overline{D}\to \overline{D'}$ and
$f^{\star}:\overline{D^{\star}}\to \overline{D'^{\star}}$ such that
$f=f^{\star}=\phi$ on $\partial D$ and $L$ depends only on $K$ and
$L_1$, where $D^{\star}=\overline{{\mathbb R}}^2\backslash
\overline{D}$ and $D'^{\star}= \overline{{\mathbb R}}^2\backslash
\overline{D'}$.
\end{Thm}

Tukia and V\"{a}is\"{a}l\"{a} \cite{TV1} dealt with the curious
phenomenon that sometimes a quasiconformal property implies the
corresponding bilipschitz property.

\begin{Thm}\label{xwm-9}$($\cite[Theorem 2.12]{TV1}$)$
Suppose that $X$ is a closed set in ${\mathbb R}^n$, $n\not=4$, and
that $f:{\mathbb R}^n \to {\mathbb R}^n$ is a $K$-QC map such that
$f|_{X}$ is $L$-bilipschitz. Then there is an $L_1$-bilipschitz map
$g:{\mathbb R}^n \to {\mathbb R}^n$ such that

$(1)$ $g|_{X}=f|_{X}$;

$(2)$ $g(D)=f(D)$ for each component $D$ of ${\mathbb R}^n\backslash
X$;

$(3)$ $L_1$ depends only on $K$, $L$ and $n$.
\end{Thm}

In \cite{Ger5}, Gehring raised the following two related problems.

\bop\label{Con2-m-1} Suppose that $D$ is a Jordan domain in
$\overline{{\mathbb R}}^2$, and that $f|_{\partial D}$ is
$M$-bilipschitz. Characterize mappings $f$  having $M'$-bilipschitz
 extension to $D$ with $M'=M'(c,M)$. \eop

 \bop\label{Con2-m-2} Suppose that $D$ is a Jordan domain in
$\overline{{\mathbb R}}^2$. For which domains $D$ does each
$M$-bilipschitz $f$ in the $\partial D$ have $M'$-bilipschitz
 extension to $D$ with $M'=M'(c,M)$? \eop

 Gehring himself discussed these two problems and got the following two
 results.

\begin{Thm}\label{xw-9}$($\cite[Theorem 2.11]{Ger5}$)$
Suppose that $D$ and $D'$ are Jordan domains  in $\overline{{\mathbb
R}}^2$ and that $\infty\in D'$ if and only if $\infty\in D$. Suppose
also that $f: D\to D'$ is a $K$-quasiconformal mapping and that $f$
extends to a homeomorphism $f: \overline{D}\to \overline{D'}$ such
that $f|_{\partial D}$ is $M$-bilipschitz. Then there exists an
$M$-bilipschitz map $g: \overline{D}\to \overline{D'}$ with
$g|_{\partial D}=f|_{\partial D}$, where $M'=M'(M, K)$.
\end{Thm}

\begin{Thm}\label{xwm-11}$($\cite[Theorem 4.9]{Ger5}$)$
Suppose that $D$ and $D'$ are Jordan domains  in $\overline{{\mathbb
R}}^2$. Then each $M$-bilipschitz $f$ in $\partial D$ has an
$M'$-bilipschitz extension $g: D\to D'$ with $g|_{\partial
D}=f|_{\partial D}$ if and only if $D$ is a $K$-quasidisk, where
$M'=M'(M, K)$ and $K=K(M)$.
\end{Thm}

We remark that Theorem \Ref{xw-9} is a partial answer to Open
Problem \ref{Con2-m-1}, and Theorem \Ref{xwm-11} is an affirmative
answer to Open Problem \ref{Con2-m-2}. In the proof of Theorem
\Ref{xw-9}, the modulus of a path family, which is an important tool
in the quasiconformal theory in ${\mathbb R}^n$, was applied. In
general, this tool is no longer applicable in the context of Banach
spaces (see \cite{Vai6-0}). A natural problem is whether Theorem
\Ref{xw-9} is true or false in Banach spaces. In fact, this problem
was raised by V\"{a}is\"{a}l\"{a} in \cite{Vai5} in the following
form.

\bop\label{Con2} Suppose that $D$ and $D'$ are bounded domains with
connected boundaries in $E$ and $E'$. Suppose also that $f: D\to D'$
is $M$-QH and that $f$ extends to a homeomorphism $f:
\overline{D}\to \overline{D'}$ such that $f|_{\partial D}$ is
$M$-bilipschitz.  Is it true that $f$  $M'$-bilipschitz with
$M'=M'(c,M)$? \eop

Our result is as follows.

\begin{thm}\label{xw-1}
Suppose that $D$ and $D'$ are bounded domains with connected
boundaries in $E$ and $E'$, respectively. Suppose also that $f: D\to
D'$ is $M$-QH and that $f$ extends to a homeomorphism $\overline{f}:
\overline{D}\to \overline{D'}$ such that $f|_{\partial D}$ is
$M$-bilipschitz. If $D'$ is a $c$-uniform domain, then $f$ is
$M'$-bilipschitz with $M'=M'(c,M)$.
\end{thm}

We see from Theorem \ref{xw-1} that the answer to Open Problem
\ref{Con2} is positive by replacing the hypothesis ``$D'$ being
bounded" in Open Problem \ref{Con2} with the one ``$D'$ being
bounded and uniform".\medskip

The organization of this paper is as follows. The proof of Theorem
\ref{xw-1} will be given in Section \ref{sec-4}. In Section
\ref{sec-2}, some preliminaries are stated including a new lemma and
its proof.

\section{Preliminaries}\label{sec-2}

The {\it quasihyperbolic length} of a rectifiable arc or a path
$\alpha$ in the norm metric in $D$ is the number (cf.
\cite{GP,Vai3}):

$$\ell_k(\alpha)=\int_{\alpha}\frac{|dz|}{d_{D}(z)}.
$$

For each pair of points $z_1$, $z_2$ in $D$, the {\it quasihyperbolic distance}
$k_D(z_1,z_2)$ between $z_1$ and $z_2$ is defined in the usual way:
$$k_D(z_1,z_2)=\inf\ell_k(\alpha),
$$
where the infimum is taken over all rectifiable arcs $\alpha$
joining $z_1$ to $z_2$ in $D$. For all $z_1$, $z_2$ in $D$, we have
(cf. \cite{Vai3})

\beq\label{eq(0000)} k_{D}(z_1, z_2)\geq
\inf\left\{\log\Big(1+\frac{\ell(\alpha)}{\min\{d_{D}(z_1), d_{D}(z_2)\}}\Big)\right\}\geq
\Big|\log \frac{d_{D}(z_2)}{d_{D}(z_1)}\Big|,\eeq
where the infimum is taken over all rectifiable curves $\alpha$ in $D$ connecting $z_1$ and $z_2$.

In \cite{Vai6},  V\"ais\"al\"a characterized uniform domains by the quasihyperbolic metric.

\begin{Thm}\label{thm0.1} {\rm (\cite[Theorem 6.16]{Vai6})}
For a domain $D$, the following are quantitatively equivalent: \bee

\item $D$ is a $c$-uniform domain;
\item $k_D(z_1,z_2)\leq c'\;
 \log\left(1+\frac{\ds{|z_1-z_2|}}{\ds\min\{d_{D}(z_1),d_{D}(z_2)\}}\right)$ for all $z_1,z_2\in D$;
\item $k_D(z_1,z_2)\leq c'_1\;
 \log\left(1+\frac{\ds{|z_1-z_2|}}{\ds\min\{d_{D}(z_1),d_{D}(z_2)\}}\right)+d$ for all $z_1,z_2\in D$.\eee
\end{Thm}

Gehring and Palka \cite{GP} introduced the quasihyperbolic metric of
a domain in $\IR^n$ and it has been recently used by many authors
 in the study of quasiconformal mappings and related questions \cite{HIMPS}.
 In the case of domains in $ {\mathbb R}^n \,,$ the equivalence
  of items (1) and (3) in Theorem D is due to Gehring and Osgood \cite{Geo} and the
  equivalence of items (2) and (3) due to Vuorinen \cite{Mvo1}. Many of the
  basic properties of this metric may be found in \cite{Geo,Vai6-0, Vai6}.

Recall that an arc $\alpha$ from $z_1$ to
$z_2$ is a {\it quasihyperbolic geodesic} if
$\ell_k(\alpha)=k_D(z_1,z_2)$. Each subarc of a quasihyperbolic
geodesic is obviously a quasihyperbolic geodesic. It is known that a
quasihyperbolic geodesic between every pair of points in $E$ exists if the
dimension of $E$ is finite, see \cite[Lemma 1]{Geo}. This is not
true in arbitrary spaces (cf. \cite[Example 2.9]{Vai4}).
In order to remedy this shortage, V\"ais\"al\"a introduced the following concepts \cite{Vai6}.

\bdefe \label{def1.4} Let $\alpha$ be an arc in $E$. The arc may be
closed, open or half open. Let $\overline{x}=(x_0,...,x_n)$, $n\geq
1$, be a finite sequence of successive points of $\alpha$. For
$h\geq 0$, we say that $\overline{x}$ is {\it $h$-coarse} if
$k_D(x_{j-1}, x_j)\geq h$ for all $1\leq j\leq n$. Let
$\Phi_k(\alpha,h)$ be the family of all $h$-coarse sequences of
$\alpha$. Set

$$s_k(\overline{x})=\sum^{n}_{j=1}k_D(x_{j-1}, x_j)$$ and
$$\ell_{k_D}(\alpha, h)=\sup \{s_k(\overline{x}): \overline{x}\in \Phi_k(\alpha,h)\}$$
with the agreement that $\ell_{k}(\alpha, h)=0$ if
$\Phi_k(\alpha,h)=\emptyset$. Then the number $\ell_{k}(\alpha, h)$
is the {\it $h$-coarse quasihyperbolic length} of $\alpha$.\edefe

In this paper, we shall use this concept in the case where $D $ is a
domain equipped with the quasihyperbolic metric $k_D$. We always use
$\ell_k(\alpha, h)$ to denote the $h$-coarse quasihyperbolic length
of $\alpha$.

\bdefe \label{def1.5} Let $D$ be a domain in $E$. An arc $\alpha\subset D$
is {\it $(\nu, h)$-solid} with $\nu\geq 1$ and $h\geq 0$ if
$$\ell_k(\alpha[x,y], h)\leq \nu\;k_D(x,y)$$ for all $x, y\in \alpha$.
A {\it $(\nu,0)$-solid arc} is said to be a {\it $\nu$-neargeodesic}, i.e.
an arc $\alpha\subset D$ is a $\nu$-neargeodesic if and only if $\ell_k(\alpha[x,y])\leq \nu\;k_D(x,y)$
for all $x, y\in \alpha$.\edefe

Obviously, a $\nu$-neargeodesic is a quasihyperbolic geodesic if and
only if $\nu=1$.

In \cite{Vai4}, V\"ais\"al\"a got the following property concerning
the existence of neargeodesics in $E$.

\begin{Thm}\label{LemA} $($\cite[Theorem 3.3]{Vai4}$)$
Let $\{z_1,\, z_2\}\subset D$ and $\nu>1$. Then there is a
$\nu$-neargeodesic in $D$ joining $z_1$ and $z_2$.
\end{Thm}

The following result due to V\"ais\"al\"a is from \cite{Vai6}.

\begin{Thm}\label{LemB2} {\rm (\cite[Theorem 4.15]{Vai6})} For domains
$D\not= E$ and $D'\not=E'$, suppose that $f: D\to D'$ is $M$-QH. If
$\gamma$ is a $c$-neargeodesic in $D$, then the arc $\gamma'$ is
$c_1$-negardesic in $D'$ with $c_1$ depending only on $c$ and
$M$.\end{Thm}
\medskip

For convenience, in the following, we always assume that $x$, $y$, $z$, $\ldots$
denote points in $D$ and $x'$, $y'$, $z'$, $\ldots$
the images in $D'$ of $x$, $y$, $z$, $\ldots$
under $f$, respectively. Also we assume that $\alpha$, $\beta$, $\gamma$, $\ldots$
denote curves in $D$ and $\alpha'$, $\beta'$, $\gamma'$, $\ldots$  the images in $D'$ of
$\alpha$, $\beta$, $\gamma$, $\ldots$
under $f$, respectively.
\medskip

\section{Bilipschitz mappings} \label{sec-4}

First we introduce the following Theorems.

\begin{Thm}\label{ThmA}$($\cite[Theorem 7.18]{Vai6}$)$
Let $D$ and $D'$ be domains in $E$ and $E'$, respectively. Suppose
that $D$ is a $c$-uniform domain and that $f:D\to D'$ is
$\varphi$-FQC $($see Section \ref{sec-2} for the definition$)$. Then
the following conditions are quantitatively equivalent:

$(1)$ $D'$ is a $c_1$-uniform domain;

$(2)$ $f$ is $\eta$-quasim\"obius.
\end{Thm}

\begin{Thm}\label{thm1.1} $($\cite[Theorem 1.1]{Ml-1}$)$
Suppose that $D$ is a $c$-uniform domain and that $f: D\to D'$ is
$(M,C)$-CQH, where $D\varsubsetneq E$ and $D'\varsubsetneq E'$. Then
the following conditions are quantitatively equivalent:

$(1)$ $D'$ is a $c_1$-uniform domain;

$(2)$  $f$ extends to a homeomorphism $\overline{f}: \overline{D}\to
\overline{D}'$ and $\overline{f}$ is $\eta$-QM rel $\partial D$.
\end{Thm}

The following theorem easily follows from Theorems \Ref{ThmA} and
\Ref{thm1.1}.

\begin{thm}\label{li-mz-1}
Suppose that $D\varsubsetneq E$ and $D'\varsubsetneq E'$, that $D$
is a $c$-uniform domain, and that $f:D\to D'$ is $\varphi$-FQC. Then
the following conditions are quantitatively equivalent:

$(1)$ $D'$ is a $c_1$-uniform domain;

$(2)$ $f$ is $\theta$-quasim\"obius;

$(3)$  $f$ extends to a homeomorphism $\overline{f}: \overline{D}\to
\overline{D}'$ and $\overline{f}$ is $\theta_1$-QM rel $\partial D$.
\end{thm}

Let us recall the following three theorems which are useful in the
proof of Theorem \ref{xw-1}.

\begin{Thm}\label{XT-0}$($\cite[Theorem 2.44]{Vai5}$)$
Suppose that $G\varsubsetneq E$ and $G'\varsubsetneq E'$ is a
$c$-uniform domain, and that $f: G\to G'$ is  $M$-QH. If $D\subset
G$  is a $c$-uniform domain, then $D'=f(D)$ is a $c'$-uniform domain
with $c'=c'(c,M)$
\end{Thm}

\begin{Thm}\label{XT-2}$($\cite[Theorem 6.19]{Vai6}$)$
Suppose that $D\varsubsetneq E$ is a $c$-uniform domain and that
$\gamma$ is a $c_1$-neargeodesic in $D$ with endpoints $z_1$ and
$z_2$. Then there is a constant $b=b(c,c_1)\geq 1$ such that

\bee
\item\label{xw-2-1} $\ds\min_{j=1,2}\ell (\gamma [z_j, z])\leq b\,d_{D}(z)
$ for all $z\in \alpha$, and

\item\label{xw-2-2} $\ell(\gamma)\leq b\,|z_{1}-z_{2}|$.
\eee
\end{Thm}

\begin{Thm}\label{lem-m-1} $($\cite[Theorem 1.2]{Yli}$)$ Suppose that $D_1$ and $D_2$ are convex domains in $E$,
where $D_1$ is bounded and $D_2$ is $c$-uniform for some $c>1$, and
that there exist $z_0\in D_1\cap D_2$ and $r>0$ such that
$\IB(z_0,r) \subset D_1\cap D_2$. If there exist constants $R_1>0$
and $c_0>1$ such that $R_1 \leq c_0r$ and $D_1\subset
\overline{\IB}(z_0,R_1)$, then $D_1\cup D_2$ is a $c'$-uniform
domain with $c'=(c+1)(2c_0+1)+c$.\end{Thm}
\medskip

\noindent {\bf Basic assumption A}\quad In this paper, we always
assume that $D$ and $D'$ are bounded domains with connected
boundaries in $E$ and $E'$, respectively, that $f: D\to D'$ is
$M$-QH, that $f$ extends to a homeomorphism $\overline{f}:
\overline{D}\to \overline{D'}$ such that $\overline{f}|_{\partial
D}$ is $M$-bilipschitz, and that $D'$ is a $c$-uniform
domain.\medskip

Before the proof  of Theorem \ref{xw-1}, we prove a series of lemmas.

\begin{lem}\label{mz-1-0}
There is a constant $M_0=M_0(M)>M$ such that for all $z_1$, $z_2\in
D$ satisfying $\dist(z_1,
\partial D)\leq \varepsilon$ and $\dist(z_2,
\partial D)\leq \varepsilon$ with sufficient small
$\varepsilon>0$, $$\frac{1}{M_0}|z_1-z_2|\leq |z'_1-z'_2|\leq
M_0|z_1-z_2|.$$
\end{lem}

\bpf Let $x_1$, $x_2\in \partial D$ be such that $|z_1-x_1|=
\frac{4}{3}\dist(z_1,
\partial D)$, $|z_2-x_2|\leq \frac{4}{3}\dist(z_2,
\partial D)$  and
$|x_1-x_2|\leq\max\{|z_1-x_1|, |z_2-x_2|\}<3|x_1-x_2|$ for
sufficient small $\varepsilon>0$. It follows from $``f$ being $M$-QH
in $D$ and homeomorphic in $\overline{D}$" that $H(x,f)\leq K$ (cf.
\cite{Vai5}) for each $x\in D$, where $K$ depends only on $M$.
Hence,

$$|z'_1-x'_1|<\frac{3}{2}K|x'_1-x'_2|\,\mbox{and}\,|z'_2-x'_2|<\frac{3}{2}K|x'_1-x'_2|.$$

If $|z_1-z_2|\leq \frac{1}{4K^2M}\max\{|z_1-x_1|, |z_2-x_2|\}$, then
for each $z\in [z_1,z_2]$, $$d_D(z)\geq
\frac{3K^2M-1}{4K^2M}\max\{|z_1-x_1|, |z_2-x_2|\},$$ and so we have
\begin{eqnarray*}
\frac{2|z'_1-z'_2|}{\min\{d_{D'}(z'_1),d_{D'}(z'_2)\}}&\leq&\log\Big(1+\frac{|z'_1-z'_2|}{\min\{d_{D'}(z'_1),d_{D'}(z'_2)\}}\Big)\\
\nonumber&\leq&
k_{D'}(z'_1,z'_2)\leq Mk_D(z_1,z_2)\\
\nonumber&\leq&M\int_{[z_1,z_2]}\frac{|dz|}{d_D(z)}
\\
\nonumber&\leq& \frac{4K^2M^2|z_1-z_2|}{(3K^2M-1)\max\{|z_1-x_1|,
|z_2-x_2|\}},\end{eqnarray*} which shows that $$|z'_1-z'_2|\leq
\frac{12K^3M^3}{3K^2M-1}|z_1-z_2|.$$

If $|z_1-z_2|> \frac{1}{4K^2M}\max\{|z_1-x_1|, |z_2-x_2|\}$, then by
the assumption $``f$ being $M$-bilipschitz in $\partial D$",
\begin{eqnarray*}|z'_1-z'_2|&\leq& |z'_1-x'_1|+|z'_2-x'_2|+|x'_1-x'_2|
\\ \nonumber&\leq& (3K+1)|x'_1-x'_2|
\\ \nonumber&\leq& (3K+1)M|x_1-x_2|\\ \nonumber&\leq& (12K+4)K^2M^2|z_1-z_2|.
\end{eqnarray*}
The same discussion as above shows that $$|z_1-z_2|\leq
(12K+4)K^2M^2|z'_1-z'_2|.$$\epf

\begin{lem}\label{mz-1-1}
There is a constant $M_1=M_1(c,M)$ such that for all $x\in D$ and
$z\in \mathbb{S}(x, d_D(x))\cap \overline{D}$ satisfying $\dist(z,
\partial D)\leq \varepsilon$ for sufficient small
$\varepsilon>0$, $$|z'-x'|\leq M_1d_{D}(x).$$
\end{lem}

\bpf Let $x_0\in \mathbb{S}(x,d_{D}(x))\cap \overline{D}$ such that
$\dist(x_0, \partial D)\leq \varepsilon$ for sufficient small
$\varepsilon>0$, and let $x_2$ be the intersection point of
$\mathbb{S}(x_0, \frac{1}{2}d_{D}(x))$ with $[x_0, x]$. Then we have
$$k_{D}(x_2,x)\leq
\log\Big(1+\frac{|x-x_2|}{d_{D}(x)-|x-x_2|}\Big)\leq \log
\frac{d_{D}(x)}{d_{D}(x_2)}=\log 2,$$ which implies that

$$\log\frac{|x'_2-x'|}{|x'_2-x'_0|}\leq k_{D'}(x'_2,x')\leq Mk_{D}(x_2,x)=M\log 2.$$ Hence
\be\label{mh-0}|x'_2-x'|\leq 2^M |x'_2-x'_0|,\ee \noindent and so
\be\label{mh-1}|x'-x'_0|\leq|x'-x'_2|+|x'_2-x'_0|
\leq(2^M+1)|x'_2-x'_0|.\ee

Let $T$ be a $2$-dimensional linear subspace of $E$ which contains
$x_0$ and $x_2$, and we use $\tau$ to denote the circle $T\cap \mathbb{S}(x_0,
\frac{1}{2}d_{D}(x))$. Take $w_1\in \tau\cap
\partial D$ such that $\tau(x_2, w_1)\subset D$ and
$\ell(\tau[x_2,w_1])\leq 2d_D(x)$. Let $x_1\in \mathbb{S}(x,
d_{D}(x))\cap \tau[x_2,w_1]\cap \overline{D}$ and denote $\tau(x_1,
w_1)$ by $\tau_1$.

\bcl\label{xt-1-1'} There must exist a $2^{32}$-uniform domain $D_1$
in $D$ and $x_3\in \partial D_1\cap \overline{D}$ satifying
$\dist(x_3,
\partial D)\leq \varepsilon$ for sufficient small $\varepsilon>0$
such that $x_0$, $x\in \overline{D}_1$ and $\frac{1}{12}d_{D}(x)\leq
|x_3-x_0|\leq \frac{11}{12}d_{D}(x)$.\ecl

If $d_{D}(x_1)=0$, then we take $D_1=\mathbb{B}(x,d_{D}(x))$ and
$x_3=x_1$.  Obviously, $ |x_3-x_0|= \frac{1}{2}d_{D}(x)$. Hence
Claim \ref{xt-1-1'} holds true in this case.

If $d_{D}(x_1)>0$, we divide the proof of Claim \ref{xt-1-1'} into
two parts.

\bca\label{li-1} $d_{D}(x_1)\leq \frac{5}{12}d_{D}(x)$.\eca

Then we take $D_1=\mathbb{B}(x, d_{D}(x))\cup \mathbb{B}(x_1,
d_{D}(x_1))$ and $x_3\in \mathbb{S}(x_1, d_{D}(x_1))\cap
\overline{D}$ such that $\dist(x_3, \partial D)\leq \varepsilon$ for
sufficient small $\varepsilon>0$. It follows from Theorem
\Ref{lem-m-1} that $D_1$ is a $29$-uniform domain and
$$\frac{1}{12}d_{D}(x)\leq |x_1-x_0|-|x_1-x_3|\leq|x_3-x_0|\leq |x_1-x_0|+|x_1-x_3|\leq\frac{11}{12}d_{D}(x),$$
from which we see that Claim \ref{xt-1-1'} is true.

\bca\label{li-1'} $d_{D}(x_1)> \frac{5}{12}d_{D}(x)$.\eca

Obviously, $d_{D}(x_1)> \frac{5}{6}|x_1-x_0|$. We let $w_2\in
\tau_1$ be the first point along the direction from $x_1$ to $w_1$
such that \be\label{mh-4}d_{D}(w_2)= \frac{5}{12}d_{D}(x).\ee

If $|w_2-x_1|\leq \frac{1}{3}d_D(x)$, then we take
$D_1=\mathbb{B}(x, d_D(x))\cup \mathbb{B}(w_2, d_D(w_2))$, and
 let $x_3\in \mathbb{S}(w_2, d_{D}(w_2))\cap \overline{D}$ such that
$\dist(x_3, \partial D)\leq \varepsilon$ for sufficient small
$\varepsilon>0$.
 Then \be\label{xv-1} \ds d_D(w_2)+d_D(x)-|w_2-x|\ge d_D(w_2)-|w_2-x_1|\ge
\frac{1}{12}d_D(x)\ee
 and \be\label{xv-2}
 \frac{1}{12}d_{D}(x)\leq |x_3-x_0|\leq
 |w_2-x_0|+|w_2-x_3|\leq\frac{11}{12}d_{D}(x).\ee

It follows from Theorem \Ref{lem-m-1} that $D_1$ is a $677$-uniform
domain, which shows that Claim \ref{xt-1-1'} is true.

If $|w_2-x_1|> \frac{1}{3}d_D(x)$, then we first prove the following
subclaim.

\bscl\label{xm-1} There exists a simply connected domain
$D_1=\bigcup \limits_{i=0}^{t}B_i$ in $D$, where $t=1$ or $2$, such
that
\begin{enumerate}
\item  $x_0$, $x\in \overline{D}_1$;
\item  For each $i\in \{0,\cdots, t\}$,  $\frac{5}{12}\,d_D(x)\leq
r_i\leq d_D(x)$;
\item If $t=2$, then $\ds |x-w_2|-r_0-r_2\ge
\frac{1}{144}d_D(x)$; and
\item  $\ds r_i+r_{i+1}-|v_i-v_{i+1}|\ge
\frac{1}{144}d_D(x)$, where $i\in\{0,1\}$ if
$t=2$ or $i=0$ if $t=1$.
\end{enumerate}
Here $B_i=\mathbb{B}(v_i, r_i)$, $v_i\in \tau[x_2,w_2]$, $v_1\not\in
B_0$ and $v_2\not\in \tau[x_2,v_1]$.\escl

To prove this subclaim, we let $y_2\in \tau_1$ be such that
$|x_1-y_2|=\frac{1}{3}d_D(x)$, and let $C_0=\mathbb{B}(x,d_D(x))$
and $C_1=\mathbb{B}(y_2, d_D(y_2))$. Since $d_D(y_2)>
\frac{5}{12}d_D(x)$, we have \be\label{xv-3}\ds
d_D(y_2)+d_D(x)-|y_2-x|\ge \frac{1}{12}d_D(x).\ee

Next, we construct a ball denoted by $C_2$:

 If
$w_2\in \overline{C}_1$, then we let $C_2=\mathbb{B}(w_2,d_D(w_2))$.

If $w_2\not\in \overline{C}_1$, then we let $y_3$ be the
intersection of $\mathbb{S}(y_2, d_D(y_2))$ with $\tau_1[y_2, w_1]$.
Since $\ell(\tau_1)\leq 2d_D(x)$ and $d_D(z)\geq \frac{5}{12}d_D(x)$
for all $z\in \tau_1(x_1,x_2)$, we have
$$|w_1-w_2|+|w_2-y_3|+|y_3-y_2|+|y_2-x_1|+|x_2-x_1|\leq
 \ell(\tau_1) \leq  2d_D(x),$$ which implies that
\be\label{xv-4}|w_2-y_3|\leq
 \frac{1}{3}d_D(x).\ee  We take $C_2=\mathbb{B}(w_2,d_D(w_2))$.
Then (\ref{xv-4}) implies \be\label{xv-5}
d_D(w_2)+d_D(x_2)-|x_2-w_2|\ge
d_D(w_2)-|w_2-y_3|\ge\frac{1}{12}d_D(x).\ee

Now we are ready to construct the needed domain $D_1$. \medskip

If $d_{D}(w_2)+d_D(x)-|w_2-x|\geq \frac{1}{48}d_{D}(x)$, then we
take $B_0=C_0$, $B_1=C_2$ and $D_1=B_0\cup B_1$ with $v_0=x$,
$v_1=w_2$, $r_0=d_D(x)$ and $r_1=d_D(w_2)$. Obviously, $D_1$
satisfies all the conditions in Subclaim \ref{xm-1}. In this case,
$t=1$.

 If
$d_{D}(w_2)+d_D(x)-|w_2-x|< \frac{1}{48}d_{D}(x)$, then we take
$B_0=\mathbb{B}(x, \frac{35}{36}d_{D}(x))$ with
$r_0=\frac{35}{36}d_{D}(x)$ and $v_0=x$, $B_1=C_1$ with
$r_1=d_D(y_2)$ and $v_1=y_2$, and $B_2=C_2$ with $r_2=d_D(w_2)$ and
$v_2=w_2$. Then Inequalities (\ref{xv-3}) and (\ref{xv-5}) show that
$D_1=\bigcup_{i=0}^{2}B_i$ satisfies all the conditions in Subclaim
\ref{xm-1}. In this case, $t=2$.

Hence the proof of Subclaim \ref{xm-1} is complete. \medskip

 It follows from a similar argument as in the proof of
\cite[Theorem $1.1$]{Ml} that

\begin{cor}\label{cor-0}
The domain $D_1$ constructed in Subclaim \ref{xm-1} is a
$2^{32}$-uniform domain.
\end{cor}

Let $x_3\in \mathbb{S}(w_2, d_{D}(w_2)\cap\overline{D}$ such that
$\dist(x_3, \partial D)\leq \varepsilon$ for sufficient small
$\varepsilon>0$. Then \be\label{mh-5}\frac{1}{12}d_{D}(x)\leq
|x_3-x_0|\leq \frac{11}{12}d_{D}(x).\ee Then the proof of Claim
\ref{xt-1-1'} easily follows from (\ref{mh-5}), Subclaim \ref{xm-1}
and Corollary \ref{cor-0}.\medskip

We come back to the proof of Lemma \ref{mz-1-1}. It follows from
(\ref{mh-5}) and Lemma \ref{mz-1-0} that \be\label{mh-6'}|x-x_3|\leq
|x-x_0|+|x_0-x_3|\leq \frac{23}{12}d_D(x)\ee and
\be\label{mh-6}\frac{1}{12M_0}d_{D}(x)\leq
\frac{1}{M_0}|x_3-x_0|\leq |x'_3-x'_0|\leq M_0|x_3-x_0|\leq
\frac{11M_0}{12}d_{D}(x).\ee Then it follows from Theorem \Ref{XT-0}
that $D'_1$ is an $M'$-uniform domain, where $M'=M'(c,M)$.  Hence we
know from Theorem \ref{li-mz-1} that $f^{-1}$ is a
$\theta$-Quasim\"obius in $\overline{D}_1$, where
$\theta=\theta(c,M)$, and so (\ref{mh-0}), (\ref{mh-1}),
(\ref{mh-5}), (\ref{mh-6'}) and (\ref{mh-6}) imply that
$$\frac{1}{23}\leq\frac{|x_3-x_0|}{|x_2-x_0|}\cdot\frac{|x_2-x|}{|x-x_3|}\leq
\theta
\Big(\frac{|x'_3-x'_0|}{|x'_2-x'_0|}\cdot\frac{|x'_2-x'|}{|x'-x'_3|}\Big)\leq
\theta\Big(\frac{M_02^{M+1}d_D(x)}{|x'-x'_3|}\Big) ,$$  which,
together with (\ref{mh-1}), shows
$$|x'-x'_0|\leq |x'-x'_3|+|x'_3-x'_0|\leq (\frac{2^{M+1}}{\theta^{-1}(\frac{1}{23})}+\frac{11M_0}{12})d_D(x)< \frac{2^{M_0+2}}{\theta^{-1}(\frac{1}{23})}d_D(x).$$
Thus the proof of Lemma \ref{mz-1-1} is complete.\epf

\begin{lem}\label{mz-1-1'}
For all $x\in D$, if $z\in \mathbb{S}(x, d_D(x))\cap \overline{ D}$
such that $\dist(x, \partial D)\leq \varepsilon$ for sufficient
small $\varepsilon>0$, then $|z'-x'|\geq \frac{1}{e^{4M_0M_1^2}}
d_{D}(x)$, where $M_1=M_1(c,M)$.
\end{lem}

\bpf  Suppose on the contrary that there exist points $x_1\in D$ and
$y_1\in \mathbb{S}(x_1,d_{D}(x_1))\cap \overline{ D}$ with
$\dist(y_1,\partial D)\leq \varepsilon$ for sufficient small $\varepsilon>0$
such that
\be\label{li-xt-0'}|x'_1-y'_1|<
\frac{1}{e^{4M_0M_1^2}}|x_1-y_1|.\ee

We take $y_2\in \mathbb{S}(y_1, d_D(x_1))\cap \overline{ D}$ such
that $\dist(y_2, \partial D)\leq \varepsilon$ for sufficient small
$\varepsilon>0$. From Lemma \ref{mz-1-0} that we know
$$|y'_1-y'_2|\geq \frac{1}{M_0}|y_1-y_2|= \frac{1}{M_0}d_D(x_1).$$

 Let $T_1$ be a
$2$-dimensional linear subspace of $E$ determined by  $x_1$, $y_1$
and $y_2$, and  $\omega$ the circle $T_1\cap \mathbb{S}(y_1,
d_{D}(x_1))$. We take $y_3\in \omega\cap \partial D$ which satisfies
$\omega(x_1, y_3)\subset D$ and $\ell(\omega[x_1,y_3])\leq
4d_D(x_1)$. Let $\omega_1=\omega(x_1, y_3)$ and $w_1$ be the first
point along the direction from $x_1$ to $y_3$ such that $$d_D(w_1)=
\frac{1}{4M_0M_1}d_D(x_1).$$

Let $v_1\in \mathbb{S}(w_1,d_{D}(w_1))\cap \overline{ D}$ such that
$\dist(w_1, \partial D)\leq \varepsilon$ for sufficient small
$\varepsilon>0$. Then it follows from Lemma \ref{mz-1-1} that
$$d_{D'}(w'_1)\leq|w'_1-v'_1|\leq M_1d_D(w_1)=\frac{1}{4M_0}d_D(x_1),$$ which,
together with  Lemmas \ref{mz-1-0} and \ref{mz-1-1} and
(\ref{li-xt-0'}),
 implies that
\begin{eqnarray*}
|x'_1-w'_1|&\geq& |y'_1-v'_1|-|x'_1-y'_1|-|v'_1-w'_1|\\ &\geq&
\frac{1}{M_0}|y_1-v_1|-\frac{1}{e^{4M_0M_1^2}}|x_1-y_1|-M_1|v_1-w_1|
\\ &\geq&
\frac{1}{M_0}(d_D(x_1)-d_D(w_1))-\frac{1}{e^{4M_0M_1^2}}|x_1-y_1|-M_1|v_1-w_1|\\
&>&\frac{1}{2M_0}d_D(x_1).\end{eqnarray*} Hence we infer from
(\ref{li-xt-0'}) that \be\label{li-xt-4}k_{D'}(x'_1, w'_1)\geq \log
\Big(1+\frac{|x'_1-w'_1|}{d_{D'}(x'_1)}\Big)> M_1^2.\ee

Since $\ell(\omega_1)\leq 4d_D(x_1)$, by the choice of $w_1$, we
have
$$k_D(x_1,w_1)\leq\int_{\omega_1[x_1,w_1]}\frac{|dx|}{d_D(x)}\leq 16M_0M_1,$$
whence $$k_{D'}(x'_1,w'_1)\leq Mk_D(x_1,w_1)\leq 16MM_0M_1,$$ which
contradicts with (\ref{li-xt-4}). The proof of Lemma \ref{mz-1-1'}
is complete.\epf

\begin{lem}\label{mz-1-2}
For  $x_1\in D$ and $x_2\in \partial D$, we have
$$|x'_1-x'_2|\leq M_2|x_1-x_2|,$$ where $M_2=2M_0+M_1$.
\end{lem}

\bpf For $x_1\in D$, we let $y_1\in \mathbb{S}(x_1,d_{D}(x_1))\cap
\overline{ D}$ such that $\dist(y_1, \partial D)\leq \varepsilon$
for sufficient small $\varepsilon>0$. Then it follows from Lemma
\ref{mz-1-1} that \be\label{mh-7}|x'_1-y'_1|\leq M_1|x_1-y_1|.\ee

For $x_2\in \partial D$, if $|y_1-x_2|\leq 2 |x_1-y_1|$, then by
Lemma \ref{mz-1-0}, we have
\begin{eqnarray*} |x'_1-x'_2|&\leq&
|x'_1-y'_1|+|y'_1-x'_2|\\ \nonumber &\leq& M_1|x_1-y_1|+M_0|y_1-x_2|\\ \nonumber &\leq&(2M_0+M_1)|x_1-y_1|\\
\nonumber &\leq& (2M_0+M_1)|x_1-x_2|.\end{eqnarray*}

If $|y_1-x_2|>2|y_1-x_1|$, then we have
$$|x_1-x_2|>|y_1-x_2|-|x_1-y_1|>\frac{1}{2}|y_1-x_2|.$$ Hence by
Lemma \ref{mz-1-0} and (\ref{mh-7}), \begin{eqnarray*}
|x'_1-x'_2|&\leq&
|x'_1-y'_1|+|y'_1-x'_2|\\ \nonumber &\leq& M_1|x_1-y_1|+M_0|y_1-x_2|\\
\nonumber &\leq& (2M_0+M_1)|x_1-x_2|,\end{eqnarray*} from which the
proof follows. \epf

\begin{lem}\label{mz-1-3}
For  $x_1\in D$ and $x_2\in \partial D$, we have
$$|x'_1-x'_2|\geq \frac{1}{M_3}|x_1-x_2|,$$ where $M_3=2M_0M_1e^{(5MM_0+8M_0)M_1^2}$.
\end{lem}

\bpf We begin with a claim.

\bcl\label{li-xt-5} For all $z\in D$, we have $d_{D'}(z')\geq
\frac{1}{e^{(5MM_0+8M_0)M_1^2}} d_{D}(z).$\ecl

To prove this claim, we let $w_2\in [z,y_1]$ be such that
$|w_2-y_1|= \frac{1}{2M_1e^{4M_0M_1^2}}d_D(z)$. It follows from
\cite{Mvo1} that
$$k_D(w_2,z)\leq\log \Big(1+\frac{|w_2-z|}{d_D(z)-|w_2-z|}\Big)<5M_0M_1^2.$$
By Lemma \ref{mz-1-1}, we have
$$|w'_2-y'_1|\leq M_1 |w_2-y_1|= \frac{1}{2e^{4M_0M_1^2}}d_D(z).$$
Hence Lemma \ref{mz-1-1'} implies $|w'_2-z'|\geq
\frac{1}{2e^{4M_0M_1^2}}d_D(z),$ whence
$$\log \frac{|w'_2-z'|}{d_{D'}(z')}\leq k_{D'}(w'_2,z')\leq
Mk_D(w_2,z)\leq 5MM_0M_1^2,$$ which shows that Claim \ref{li-xt-5}
is true.
\medskip

Now we are ready to finish the proof of Lemma \ref{mz-1-3}. For $x_1\in D$ and
$x_2\in
\partial D$, if $|x_1-x_2|\leq 2M_0M_1d_D(x_1)$, then by Claim \ref{li-xt-5},
$$|x'_1-x'_2|\geq d_{D'}(x'_1)\geq\frac{1}{e^{(5MM_0+8M_0)M_1^2}}
d_{D}(x_1)\geq \frac{1}{2M_0M_1e^{(5MM_0+8M_0)M_1^2}}|x_1-x_2|.$$

If $|x_1-x_2|>  2M_0M_1d_D(x_1)$, then we take $w_3\in
\mathbb{S}(x_1, d_D(x_1))\cap \overline{ D}$ such that $\dist(w_3,
\partial D)\leq \varepsilon$ for sufficient small $\varepsilon>0$,
and so
$$|w_3-x_2|\geq |x_1-x_2|-|x_1-w_3|\geq \Big(1-\frac{1}{
 2M_0M_1}\Big)|x_1-x_2|$$
 and $$|w_3-x_2|\geq |x_1-x_2|-|x_1-w_3|\geq (2M_0M_1-1)|x_1-w_3|,$$
whence Lemmas \ref{mz-1-0} and \ref{mz-1-1} implies
 \begin{eqnarray*} |x'_1-x'_2|&\geq&|w'_3-x'_2|-|x'_1-w'_3|\\
 &\geq& \frac{1}{M_0}|w_3-x_2|-M_1|x_1-w_3|\\
 &\geq& \Big(\frac{1}{M_0}-\frac{M_1}{2M_0M_1-1}\Big)|w_3-x_2|
 \\
 &\geq& \frac{1}{3M_0}|x_1-x_2|,\end{eqnarray*}
from which the proof is complete. \epf

By the above lemmas, we get the following result.

\begin{lem}\label{mz-1-4}
$D$ is a $c_1$-uniform domain, where $c_1=c_1(c,M)$.
\end{lem}

\bpf We first prove that $f^{-1}$ is  $\theta_1$-Quasim\"obius rel
$\partial D'$, where $\theta_1(t)=(M_2M_3)^2t$, $M_2$ and $M_3$ are
the same as in Lemma \ref{mz-1-2} and Lemma \ref{mz-1-3},
respectively. By definition, it is necessary to prove that for
$x'_1$, $x'_2$, $x'_3$, $x'_4\in \overline{D'}$, \beq\label{Fri-1}
\frac{|x_4-x_1|}{|x_4-x_2|}\cdot\frac{|x_2-x_3|}{|x_1-x_3|}
  &\leq& (M_2M_3)^2\frac{|x'_4-x'_1|}{|x'_4-x'_2|}\cdot\frac{|x'_2-x'_3|}{|x'_1-x'_3|},\eeq where $x_1,x_2\in
\partial D'$. Obviously, to prove Inequality \ref{Fri-1}, we only
need to consider the following three cases.

\bca\label{mz-li-xt-1-1''}$x'_1,x'_2,x'_3, x'_4\in
\partial D'$.\eca

Since $f$ is $M$-bilipschitz in $\partial D$, we have
\begin{eqnarray*}\frac{|x_4-x_1|}{|x_4-x_2|}\cdot\frac{|x_2-x_3|}{|x_1-x_3|}
  &\leq& M^4\frac{|x'_4-x'_1|}{|x'_4-x'_2|}\cdot\frac{|x'_2-x'_3|}{|x'_1-x'_3|}.\end{eqnarray*}

\bca\label{mz-li-xt-1-1'}$x'_1$, $x'_2$, $x'_3\in \partial D'$ and
$x'_4\in D'$.\eca

It follows from Lemmas \ref{mz-1-2} and \ref{mz-1-3} that
\begin{eqnarray*}\frac{|x_4-x_1|}{|x_4-x_2|}\cdot\frac{|x_2-x_3|}{|x_1-x_3|}
 &\leq& \frac{M_2M_3|x'_4-x'_1|}{|x'_4-x'_2|}\cdot\frac{M^2|x'_2-x'_3|}{|x'_1-x'_3|}
  \\ &\leq& M^2M_2M_3\frac{|x'_4-x'_1|}{|x'_4-x'_2|}\cdot\frac{|x'_2-x'_3|}{|x'_1-x'_3|}.\end{eqnarray*}

\bca\label{mz-li-xt-2} $x'_1$, $x'_2\in \partial D'$ and $x'_3,
x'_4\in D'$.\eca

We obtain from Lemmas \ref{mz-1-2} and \ref{mz-1-3} that
\begin{eqnarray*}\frac{|x_4-x_1|}{|x_4-x_2|}\cdot\frac{|x_2-x_3|}{|x_1-x_3|}
 &\leq& \frac{M_2M_3|x'_4-x'_1|}{|x'_4-x'_2|}\cdot\frac{M_2M_3|x'_2-x'_3|}{|x'_1-x'_3|}
  \\ &\leq& (M_2M_3)^2\frac{|x'_4-x'_1|}{|x'_4-x'_2|}\cdot\frac{|x'_2-x'_3|}{|x'_1-x'_3|}.\end{eqnarray*}

The combination of Cases \ref{mz-li-xt-1-1''} $\sim$
\ref{mz-li-xt-2} shows that Inequality \eqref{Fri-1} holds, which
implies that $f^{-1}$ is a $\theta_1$-Quasim\"obius rel $\partial
D'$. Hence Theorem \ref{li-mz-1} shows that $D$ is a $c_1$-uniform
domain, where $c_1$ depends only on $c$ and $M$. \epf

\subsection{ The proof of Theorem \ref{xw-1}} For any $z_1$,
$z_2\in \overline{D}$, it suffices to prove that
\beq\label{Sat-2} \frac{1}{M'}|z_1-z_2|\leq |z'_1-z'_2|\leq M'|z_1-z_2|,\eeq where
$M'$ depends only on $c$ and $M$.

It follows from the hypothesis ``$f$ being $M$-bilipschitz in $\partial D$", Lemmas
\ref{mz-1-2} and \ref{mz-1-3} that we only need to consider the case
$z_1$, $z_2\in D$.

If $|z_1-z_2|\leq \frac{1}{2}\max\{d_D(z_1), d_D(z_2)\}$, then
$$k_D(z_1,z_2)\leq \int_{[z_1,z_2]}\frac{|dx|}{d_D(x)}\leq \frac{2|z_1-z_2|}{\max\{d_D(z_1), d_D(z_2)\}}\leq 1,$$
which shows that $$\log\Big(1+\frac{|z'_1-z'_2|}{\min\{d_{D'}(z'_1),
d_{D'}(z'_2)\}}\Big)\leq k_{D'}(z'_1,z'_2)\leq Mk_D(z_1,z_2)\leq
M,$$
and so \beq\label{mh-12}\frac{|z'_1-z'_2|}{e^M\min\{d_{D'}(z'_1),
d_{D'}(z'_2)\}}&\leq&
\log\Big(1+\frac{|z'_1-z'_2|}{\min\{d_{D'}(z'_1),
d_{D'}(z'_2)\}}\Big)\\ \nonumber&\leq&
\frac{2M|z_1-z_2|}{\max\{d_D(z_1), d_D(z_2)\}}.\eeq We see from
Lemma \ref{mz-1-1} that $$\min\{d_{D'}(z'_1), d_{D'}(z'_2)\}\leq
M_1\max\{d_D(z_1), d_D(z_2)\}.$$
Then (\ref{mh-12}) implies that \be\label{mh-14}|z'_1-z'_2|\leq
2MM_1e^M|z_1-z_2|.\ee

For the other case $|z_1-z_2|> \frac{1}{2}\max\{d_D(z_1),
d_D(z_2)\}$, we let $\beta$ be a $2$-neargeodesic joining $z_1$ and
$z_2$ in $D$. It follows from Theorem \Ref{LemB2} that $\beta'$ is a
$c_2$-neargeodesic, where $c_2$ depends only on $M$. Let
$z'\in\beta'$ such that
$$|z'_1-z'|=\frac{1}{2}|z'_1-z'_2|.$$ Then we know from
$|z'_2-z'|\geq \frac{1}{2}|z'_1-z'_2|$ and Theorem \Ref{XT-2} that
\beq\label{li-x-1}|z'_1-z'_2|&\leq& 2\min\{|z'_1-z'|, |z'_2-z'|\}\\
\nonumber&\leq& 2\min\{\diam(z'_1,z'), \diam(z'_2,z')\}\\
\nonumber&\leq& 2\mu d_{D'}(z'),\eeq where $\mu$ depends only on $c$
and $M$.

We claim that \be\label{mh-13}d_D(z)\leq 3\ell(\beta).\ee Otherwise,
$$\max\{d_D(z_1),
d_D(z_2)\}\geq d_D(z)-\max\{|z_1-z|, |z_2-z|\}> 2\ell(\beta)\geq
2|z_1-z_2|.$$ This is the desired contradiction.

By Theorem \Ref{XT-2} and Lemma \ref{mz-1-4}, we have $$d_{D}(z)\leq
3 \ell(\beta)\leq 3b|z_1-z_2|,$$ where $b=b(c_1)$. Hence Lemma
\ref{mz-1-1} and (\ref{li-x-1}) show that
\be\label{mh-15}|z'_1-z'_2|\leq 2\mu d_{D'}(z')\leq
6bM_1\mu|z_1-z_2|.\ee

By Lemma \ref{mz-1-4}, we see that $D$ is a $c_1$-uniform domain.
Hence a similar argument as in the proofs of Inequalities
(\ref{mh-14}) and (\ref{mh-15}) yields that \beq
\label{Sat-1}|z_1-z_2|\leq M_4|z'_1-z'_2|,\eeq where $M_4=M_4(c,M)$.

Obviously, the inequalities \eqref{mh-14}, \eqref{mh-15} and
\eqref{Sat-1} show that \eqref{Sat-2} holds, and thus the proof of
the theorem is complete. \qed

\end{document}